\documentclass[12pt]{elsarticle}
\usepackage{amsmath}
\usepackage{helvet}
\usepackage{amssymb}
\usepackage[margin=1 in]{geometry}
\usepackage{amsthm}
\usepackage{lipsum}
\usepackage{chngcntr}
\PassOptionsToPackage{square, numbers}{natbib}
\usepackage[square,numbers]{natbib}
\usepackage{relsize}
\usepackage{titlesec}
\newtheorem{thm}{Theorem}
\newtheorem{lem}{Lemma}
\newtheorem{definition}{Definition}
\newtheorem{cor}{Corollary}
\begin{document}
 \title{\textbf{Infinitely divisible nonnegative matrices,\\ $M$-matrices,\\ and the embedding problem for finite state stationary Markov Chains}}
 \date{06 August 2015}
 \author{Alexander Van-Brunt}
 \ead{alexvb@kurims.kyoto-u.ac.jp}
 
 \address{ \mbox{Research Institute for Mathematical Sciences}\\
Kyoto University; \\
Kyoto 606-8502 JAPAN FAX: +81-75-753-7272}


 \begin{abstract}
 This paper explicitly details the relation between $M$-matrices, nonnegative roots of nonnegative matrices, and the embedding problem for finite-state stationary Markov chains. The set of nonsingular nonnegative matrices with arbitrary nonnegative roots is shown to be the closure of the set of matrices with matrix roots in $\mathcal{IM}$. The methods presented here employ nothing beyond basic matrix analysis, however it answers a question regarding $M$-matrices posed over 30 years ago and as an application, a new characterization of the set of all embeddable stochastic matrices is obtained as a corollary. \\\\
 \textit{AMS Classification}:   15B51,  15B48
 \end{abstract}

 \begin{keyword} $M$-Matrix \sep Markov chains \sep Embedding problem \sep Nonnegative matrix \sep Matrix roots \end{keyword}
 
 \maketitle
 \newpage

 \section*{Introduction} 
A $Z$-matrix is a matrix that has nonpositive off-diagonal elements. An $M$-matrix is defined as a $Z$-matrix that has a non-negative inverse, or alternatively, as a $Z$-matrix of the form $\alpha I-K$.  Here, $K$ is a nonnegative matrix and $\alpha > \rho(K)$, where $\rho$ is the spectral radius of $K$. In fact there are many characterizations of $M$-matrices. A summary of these characterizations can be found in \cite{Plemmons1977} and \cite{MatrixAnalysis2}. The `Inverse $M$-matrix problem' concerns the conditions under which a nonnegative matrix is the inverse of an $M$-matrix. The set of such nonnegative matrices is denoted by $\mathcal{IM}$. Two key surveys of this problem are given by C.Johnson in \cite{Inverse$M$-matrices1}, and more recently C.Johnson and R.Smith in \cite{InverseM-matrices2}.
More then 30 years ago, the question was raised by C.Johnson in \cite{Inverse$M$-matrices1} asking for which nonnegative  matrices $B$, does there exist a sequence of nonnegative matrices $\lbrace K_{n} \rbrace_{n=1}^{\infty}$ such that  

\begin{equation}
 (K_{n})^{n}=B.
\end{equation}

\noindent Informally put; which nonnegative matrices possess nonnegative matrix roots of arbitrary order.  Indeed the question in \cite{Inverse$M$-matrices1} asked specifically if a nonsingular, non-negative matrix that has arbitrary, nonnegative roots,  also has roots which are in $\mathcal{IM}$.  We shall see this statement is correct, modulo some further conditions.  \\

 The question in \cite{Inverse$M$-matrices1} is connected with the embedding problem for Markov chains. The latter has been a long standing problem in linear algebra and probability theory since it was first considered by Elfving \cite{Elfving}. The precise formulation of this problem will be given later; however, a connection was made by Kingman in \cite{Kingman1962}, who showed that a Markov chain is embeddable if and only if its stochastic matrix is nonsingular and has arbitrary stochastic matrix roots.  \\
 
 Extensive work was also done on a analogous problem viz the characterization of the class of nonnegative definite matrices having the property that every positive fractional Hadamard power is also nonnegative definite. This was pioneered in \cite{RAHorn1969}. In the context of nonnegative definiteness, such matrices are called infinitely divisible. Following this terminology we make the following definition.\\
 
 \begin{definition}
  A nonnegative matrix $B$ is said to be infinitely divisible if and only if there exists a sequence of nonnegative matrices $\lbrace K_{n} \rbrace_{n=1}^{\infty}$ such that

\begin{equation} \label{ID}
 (K_{n})^{n}=B.
\end{equation}

\noindent If in addition $\det (B) > 0$ we say that $B$ is strongly infinitely divisible. 
\end{definition}

 In this paper we develop a theory for these classes of matrices and answer the question in \cite{Inverse$M$-matrices1}. Although the results will be of interest in other fields,  the embedding problem for finite state stationary Markov chains is the primary application intended in this paper. \\
\indent The paper is organized as follows: the first section is a statement of the main results; the second and third section detail the proofs and framework; and the forth is dedicated to the embedding problem. The reader who is only interested in the results for the embedding problem may thus proceed directly to section $4$; the reader interested primarily in linear algebra may omit section $4$ altogether. 
 
 \bigskip 
 
\section{Main results}

  The first result is a characterization of the set of strongly infinitely divisible matrices in terms of the exponential map.\\

  \begin{thm} \label{Z-matrix theorem}
  An $N \times N$ nonnegative matrix $B$ is strongly infinitely divisible if and only if  there exists a $Z$ matrix, $Q$ such that $e^{-Q}=B$.
\end{thm}

 \noindent This result will be used to answer the question in \cite{Inverse$M$-matrices1}. \\
  
\begin{thm} \label{M-matrix theorem}
 The set of infinitely divisible nonnegative matrices contains the closure of the set

\begin{equation}
 \lbrace B \: : \: B=K^{n}, \: n \in \mathbb{N}, \: K \:  \in \mathcal{IM}\rbrace. 
\end{equation}

\noindent Furthermore, if $B$ is nonsingular, it is infinitely divisible if and only if it belongs in the closure of this set. 
 \end{thm}

 \noindent We also prove the following result that relates the strongly infinite divisibility of the matrix to that of its submatrices.\\
 
 \begin{thm} \label{dichotomy theorem}
Let $B$ be a strongly infinitely divisible nonnegative matrix. Then following dichotomy holds:\\
\begin{enumerate}
\item If $B$ is irreducible, it is strictly positive \\

\item If $B$ is reducible, then there exists a permutation matrix $L$ such that $B=L U L^{T}$ for some infinitely divisible upper block triangular matrix $U$. Furthermore:\\

\begin{enumerate}

\item  All the matrices on the diagonal blocks of $U$ are strictly positive. \\
\item  If $M$ is the number of blocks in $U$, then for every $n \leq M$ the square submatrix, $U^{(n)}$, obtained by deleting the first $n$ blocks from the top rows and left columns, is strongly infinitely divisible.

\end{enumerate}

\end{enumerate}

\end{thm}

 We prove Theorems \ref{Z-matrix theorem} and \ref{M-matrix theorem} in section 2 and Theorem \ref{dichotomy theorem} in the subsequent one, following a number of other algebraic properties.  The main results are given above, but there are several other results contained in the following sections that are also of their own interest. 
  
\section{Connection to $Z$-matrices and Inverse $M$-matrices}

The initial part of this analysis is along very similar lines to the work by Kingman in \cite{Kingman1962}; however it is not exclusive to stochastic matrices. We begin with a proof a Theorem \ref{Z-matrix theorem}.\\

\begin{proof}[Proof of Theorem \ref{Z-matrix theorem}] 
A proof of the direct implication can be found in \cite{MatrixFunction}. It is presented here for the readers convenience. Suppose that $B=e^{-Q}$ for some $Z$-matrix $Q$. Then an $n$th root of $B$ is $e^{-Q / n}$, which is again the exponential of the negative of a $Z$-matrix. Non
negativity follows by taking a sufficiently large $\theta$ so that $-Q / n +\theta I$ is non-negative and then writing 

\begin{equation} 
e^{-Q / n} = e^{-\theta} e^{-Q / n +\theta I}, 
\end{equation}

\noindent where $I$ denotes the identity matrix.\\

\indent Conversely, assume that $B$ is a $N \times N$ strongly infinitely divisible matrix.  We first show that the sequence $B^{\frac{1}{n}}$ contains a subsequence converging to the $N \times N$ identity matrix. Let $M$ be an integer that is divisible for every integer $k$ less then or equal to $N$. Define the sequence $R_{n}= B^{\frac{1}{Mn}} \geq 0$. The relation $R^{Mn}_{n} =B$ implies that $R_{n}$ is bounded and so will have a convergent subsequence, say $R_{n_{j}}$ with a limit $R$. By the Perron Frobenius theorem, $B$ has a strictly positive eigenvalue $\widetilde{\lambda}$  which is the spectral radius of $B$. Elementary considerations tell us that for $R_{n}$ to be a real, let alone a nonnegative root of $B$, we must take the real root of $\widetilde{\lambda}$. By the same reasoning $\widetilde{\lambda}^{1/Mn}$ is the spectral radius of $R_{n}$. Hence the spectral radius of $R$ is $1$. A similar argument shows that the determinant of $R$ is on the unit circle, indeed;

\begin{equation} 
|\det (R_{n}) | = | \det ( B ) |^{\frac{1}{Mn}} \rightarrow 1, \:\:\:\:\:\: \text{as} \: n \rightarrow \infty.
\end{equation}

 \noindent Therefore every eigenvalue of $R$ must be on the unit circle (at this point nonsingularity is essential). By assumption, $R$ is non-negative so by the Perron Frobenius theorem if $R$ is irreducible then every eigenvalue is a root of unity for some $k \leq N$. If $R$ is not irreducible then we may decompose $R$ into the form $PWP^{T}$ where $P$ is a permutation matrix and $W$ is a block upper triangular matrix in which each diagonal block is an irreducible nonnegative matrix. As the spectrum of $R$ is the union of the spectra of the diagonal blocks in $W$, every eigenvalue is a $k$th root for some integer less then $n \leq N$.  Therefore in all cases, $R^{M}=I$. For a review on the Perron Frobenius theorem and irreducible matrices, we direct the reader to \cite{MatrixAnalysis1}. We thus conclude that there is a subsequence $\lbrace n_{j} \rbrace$ such that

\begin{equation}
 \lim_{j \rightarrow \infty} B^{\frac{1}{M n_{j}}} = I.
  \end{equation}
  
  We can now estimate the decay rate of the diagonal elements of $B^{1/Mn}$ using the inequality
  
  \begin{equation} \label{decayineq}
  |\text{Tr} (B^{1/Mn})|=\Big |\sum_{i=1}^{N} \lambda_{i}^{(1/Mn)} \Big | \leq N \widetilde{\lambda}^{1/Mn} = N+\mathcal{O}(n^{-1}). 
  \end{equation}
  
  \noindent Set $k=M n_{j}$, and $B^{\frac{1}{k}}=I+A_{k}$ where $A_{k}$ is some sequence of matrices with nonnegative off diagonal elements converging to $0$ (here and henceforth convergence is in any matrix norm). From inequality \eqref{decayineq} we know that the diagonal elements on the matrix $A_{k}$ decay as $\mathcal{O}(k^{-1})$ as $k \rightarrow \infty$. We also have that

  \begin{equation} \label{Stable}
  \big ( I+A_{k} \big )^{k} = \sum_{m=1}^{k} \binom{k}{m} \big [A_{k} \big ]^{m}=B. 
  \end{equation}
  
  Equation \eqref{Stable} along with the fact that all off diagonal elements of $A_{k}$ are nonnegative and $\binom{k}{m}$ is $\mathcal{O}(k^{m})$ implies that all elements of $A_{k}$ decay as $\mathcal{O}(k^{-1})$. Indeed, if an off diagonal element were to decay slower, then there would be no negative term in the diagonal to match this slower decay rate and maintain the relationship in (\ref{Stable}). We thus have that the sequence

\begin{equation} 
C_{k}:=k(B^{\frac{1}{k}}-I)=kA_{k}
\end{equation}

\noindent is bounded. Hence we find another subsequence of $C_{k}$, denoted $C_{q_{k}}$ that is convergent to some limit $C$. 
We show that $e^{C}=B$. Let

\begin{equation}
q_{k}(B^{\frac{1}{q_{k}}}-I)=C+\epsilon_{q_{k}}
\end{equation}

\noindent For some matrices $\epsilon_{q_{k}}$, that converge to $0$. Rearranging, we have for each $q_{k}$

\begin{equation}
B=\Big ( \frac{C}{q_{k}}+\frac{\epsilon_{q_{k}}}{q_{k}} +I \Big )^{q_{k}}
\end{equation}

and the relation follows from standard estimates and the binomial theorem.  \\\\
\indent Furthermore, because each $A_{q_{k}}$ has nonnegative off diagonal elements, $A$ has nonnegative off diagonal elements. Taking $Q=-A$ we complete the proof. \end{proof}

\bigskip
An interesting implication of Theorem \ref{Z-matrix theorem} is that the nonnegative roots of infinitely divisible matrices cannot be scattered: they must belong to the same branch of roots. The primary difficultly in dealing with singular infinitely divisible matices is that nothing similar to Theorem \ref{Z-matrix theorem} seems to apply. For example, the zero matrix. \\
\indent In general there is no uniqueness of the $Z$ matrix. In fact, there may be an uncountable family of $Z$ matrices associated an infinitely divisible matrix. In the context of the embedding problem, such an example is provided in \cite{Singer1976}. \bigskip

 A result central to this paper is Theorem $11$  in \cite{Inverse$M$-matrices1} which states that the primary $n$th root of an $M$-matrix is also an $M$-matrix. This yields in the following (also noted in \cite{Inverse$M$-matrices1}): \\
 
 \begin{thm}
  Let $B$ be a nonnegative matrix such that $B=K^{n}$ for some $K \in \mathcal{IM}$ then $B$ is is strongly infinitely divisible.
 \end{thm}
 
To prove theorem $2$ we also need the following lemma:\\

\begin{lem}  \label{closedlemma}
The set of infinitely divisible matrices is closed. 
\end{lem}

\begin{proof} 
Let $\lbrace B_{n} \rbrace$ be a sequence of infinitely divisible matrices  converging to $B$. Then, for any given $m \in \mathbb{N}$, we can consider the sequence $B_{n}^{1/m} \geq 0$. This sequence contains a convergent subsequence with some limit $\widetilde{B}$. By continuity $\widetilde{B}$ is nonnegative and satisfies $\widetilde{B}^{m}=B$. 
\end{proof} 
\bigskip

\begin{proof}[Proof of Theorem \ref{M-matrix theorem}]

Assume there is a sequence $\lbrace B_{n} \rbrace_{n=1}^{\infty}$ converging to $B$ such that each $B_{n}$ is the power of some inverse $M$-matrix. This implies each $B_{n}$ is infinitely divisible and so by Lemma \ref{closedlemma}, we have that $B$ is also infinitely divisible. This completes the first statement of Theorem \ref{M-matrix theorem}.\\\

Conversely suppose that $B$ is strongly infinitely divisible and that the off diagonal entries of the associated $Z$-matrix, $Q$, are strictly negative. Then,

\begin{equation} B^{-1}=e^{Q}. \end{equation}

\noindent Taking $n$th roots for $n$ sufficently large we see

\begin{equation} 
B^{-\frac{1}{n}}=I+\frac{Q}{n}+O \big ( n^{-2} \big ).
\end{equation}

\noindent is an $Z$-matrix and specifically because it's inverse is positive, it is also an $M$-matrix. It is clear that such matrices; are dense in the set of infinitely divisible matrices, hence, we may apply Lemma \ref{closedlemma} to prove the result. \end{proof}

\bigskip
The above result also indicates that the infinitely divisible matrices, whose associated $Z$-matrix  have no off-diagonal zeros, always have roots in $\mathcal{IM}$. Once one violates this condition, it is easy to construct strongly infinitely divisible matrices that do not have roots in $\mathcal{IM}$. \\

How might powers of inverse $M$-matrices be characterized? If $B$ is indeed the power of an inverse $M$-matrix. Then using the series expansion of $(1+x)^{-n}$ we see that $B$ must be of the form
\begin{equation} \sum_{k=0}^{\infty} \binom{n+k-1}{k} P^{k}, \end{equation}
 for some nonnegative matrix $P$ where, in order for the series to converge, $\rho(P)<1$. Conversely if $B$ it is of the above form for some $n >0$, then it must be power of an inverse $M$-matrix.\\
 
 We will also require the following lemmata.\\
 
 \begin{lem} \label{infdivdense}
The set of matrices with distinct eigenvalues is dense in the set of strongly infinitely divisible nonnegative matrices.
\end{lem}

\begin{proof}
Consider the set of all $Z$-matrices. This set is a convex and satisfies the condition of Corollary $2$ in \cite{Denseeigenvalues}, which implies that the set of matrices with distinct eignvalues is dense in this set. The density of the eigenvalues in the set of infinitely divisible matrices is now simply a consequence of the continuity of the exponential map on matrices.\end{proof}  
\bigskip

We also recall the following fact. See \cite{MatrixAnalysis2}.\\

\begin{lem}
Every $Z$-matrix has a nonnegative eigenvector.  
\end{lem}

\bigskip
 We conclude this section by providing a bound on the eigenvalues of $Z$-matrices generating a infinitely divisible matrix through the exponential map. Thus if one is checking for the existence of said $Z$-matrices for a nonnegative matrix with distinct eigenvalues, one need  only check a finite number of them.\\

 \begin{thm} \label{bound on eigenvalues}
 Let $B=e^{-Q}$ be a strongly infinitely divisible nonnegative matrix. Let $v$ be a nonnegative eigenvector of $Q$ with associated eigenvalue $\widetilde{\lambda}$. 
 Then every eigenvalue $ \lambda $ of $Q$ satisfies

\begin{equation} | \textnormal{Im} (\lambda)| \leq |\widetilde{\lambda} (n-1) - \log( \det (B) ) |. \end{equation}

   \end{thm} 
 \bigskip

 \begin{proof}
 Assume that $B=e^{-Q}$ for some $Z$ matrix $Q$, with entries $Q_{ij}$. Then

 \begin{equation} (\widetilde{\lambda} +Q_{ii})v_{i}  = - \sum_{j \neq i} Q_{ij} v_{j} \geq 0. \end{equation}

 \noindent In particular we have that $ -Q_{ii} \leq  \widetilde{\lambda}$, which provides an upper bound on the diagonal elements. Combining this with

 \begin{equation} \prod e^{-Q_{jj}} = \det (B), \end{equation}

 \noindent we can deduce that $-Q_{ii} \geq -\widetilde{\lambda} n +\log( \det (B) ) $, and thus we conclude that the matrix

 \begin{equation} -Q+\big (\widetilde{\lambda} (n-1) - \log( \det (B) )  \big ) I \end{equation}

 \noindent is nonnegative with nonnegative eigenvalue $\widetilde{\lambda} (n-1) - \log( \det (B) )$. By the Perron-Frobenius theorem, this is its spectral radius; therefore,

 \begin{equation} | \textnormal{Im} (\lambda)| \leq |\widetilde{\lambda} (n-1) - \log( \det (B) ) |. \end{equation}
 \end{proof}

 The above results and proofs are similar in spirit to those presented in \cite{Israel2001} and we shall see that when dealing with a stochastic matrix, the bound simplifies considerably.

 \bigskip
 \bigskip
 

\section{Algebraic properties.} 


 Given the known invariant zero patterns for $M$-matrices, it is natural to ask similar questions regarding infinitely divisible matrices. The following result is one such immediate observation:\\

\begin{lem}  \label{zero pattern invariancelem}
If, for an strongly infinitely divisible matrix $B$, $B_{ij}=0$ for $i \neq j$, then on any $Z$ matrix $Q$ such that $B=e^{-Q}$, we have $(Q^{n})_{ij}=0$ for every $n \in \mathbb{N}$. 
\end{lem}

\begin{proof}
 Let $B=e^{-Q}$, where $Q$ is a $Z$ matrix and write $Q=\theta I +(Q-\theta I)$, where $\theta$ is chosen sufficiently large so that $(Q-\theta I)$ has negative diagonal elements. Then

\begin{equation} 
B=e^{-Q}=e^{(\theta I -Q)-\theta I}=e^{-\theta }  e^{\theta I - Q}.
\end{equation}

Now $e^{-\theta }  > 0$,  so for $B_{ij}$ to be zero, we must have that

\begin{equation}
 (e^{\theta I - Q})_{ij} =0.
  \end{equation}

\noindent But this a positive power series of nonnegative matrices; therefore, for $B_{ij}$ to be $0$

\begin{equation} \label{zeropow}
[(\theta I-Q)^{n}]_{ij}=0.
\end{equation} 

\noindent for every $n$. In particular \eqref{zeropow} implies that $(Q^{n})_{ij}=0$ for each $n$ (as the term we are considering is off diagonal). \end{proof} 

Lemma \ref{zero pattern invariancelem} yields the following corollaries: \\

\begin{cor} \label{zero pattern invariance}
For a strongly infinitely divisible nonnegative matrix $B$, $B_{ij}=0$ ($i \neq j$)  implies that $(B^{t})_{ij}=\big ( e^{-tQ} \big )_{ij}=0$ for all $t \in \mathbb{R}$. 
\end{cor}

Along similar methods, one can prove other similar statements, for instance taking $\theta$ sufficiently large in (\ref{zeropow}) we can deduce. 
\\

 \begin{cor} A strongly infinitely divisible matrix has strictly positive diagonal elements. \end{cor}

In order to `test' if a matrix is infinitely divisible, it is of interest to know what operations infinitely divisible nonnegative matrices are closed under. More generally, how we might alter an infinitely divisible matrix such that it remains infinitely divisible? Infinitely divisible matrices are not closed under addition. Even in the $2 \times 2$ case consider the example where

\begin{equation}
 A=\begin{pmatrix} 2 & \frac{6}{5} \\ 3  & 2 \end{pmatrix},
 \end{equation}
 
\noindent then $A$ and $A^{T}$ are both infinitely divisible but $\det (A+A^{T}) < 0$. This example also shows that the set of infinitely divisible matrices is not convex.\\

We can also show that the product of two infinitely divisible stochastic matrices need not be infinitely divisible. As a specific counterexample, we show the product of two embeddable stochastic matrices need not be embeddable. 
 Consider the two intensity matrices
 
 \begin{equation} 
 Z_{1}=\begin{pmatrix} 
 -2 & 1 & 1 \\ 0 & -1 & 1\\ 0& 0&0 \end{pmatrix},
  \end{equation}

\begin{equation} 
Z_{2}=\begin{pmatrix}
 -\frac{1}{2} & \frac{1}{12} & \frac{5}{12} \\ 0 & -3 & 3\\ 0& 0&0 \end{pmatrix}.
 \end{equation}
 
 Then the associated stochastic matrices are:

\begin{equation} \label{countexam1}
E_{1}:=\exp (Z_{1}) \approx \begin{pmatrix}
  0.135 & 0.233 & 0.632 \\ 0 & 0.368 & 0.632\\ 0& 0&1
  \end{pmatrix},
  \end{equation}
  
 and

\begin{equation} \label{countexam2}
E_{2}=\exp (Z_{2}) \approx \begin{pmatrix}
  0.607 & 0.018 & 0.375 \\ 0 & 0.050 & 0.950\\ 0& 0&1 \end{pmatrix},
  \end{equation}
  
  The principle branch of the logarithm of $E_{2} E_{1}$ yields a matrix that negative off diagonal entries. Furthermore,  this matrix has distinct positive eigenvalues, so that the only possible intensity matrix that can generate $E_{2} E_{1}$ is the principal branch of the logarithm.  Therefore the matrix $E_{2} E_{1}$ is not embeddable.\\
  
Curiously though, $E_{1}$ and $E_{2}$ are elements of  $\mathcal{IM}$ and hence infinitely divisible. Furthermore it is easily shown that
$E_{1}E_{2} \in \mathcal{IM}$, so that $E_{1}E_{2}$ is embeddable. In summary,\\

\textit{Let $A,B$ be two infinitely divisible matrices. Then:
\begin{itemize} 
\item $AB$ need not be infinitely divisible.
\item The product of $M$-matrices need not be the power of an $M$-matrix.
\item If $AB$ is infinitely divisible, then $BA$ need not be.
\item If $AB$ is the power of an $M$-matrix, $BA$ need not be.
 \end{itemize} }
 
\indent There is however an important class of embeddable matrices for which the product of them is again embeddable.\\

\begin{thm}
 Let $A$ and $B$ be two commuting infinitely divisible matrices. Then $AB=BA$ is infinitely divisible. 
 \end{thm}
 
\begin{proof}
We can without loss of generality suppose that $A,B$ have distinct eigenvalues, and then use Lemma \ref{closedlemma}. In this case, any matrix function of $A$ or $B$ is primary, and therefore is a polynomial of $A$ or $B$ respectively (see \cite{$M$-MatrixFunctions} for more details regarding matrix functions). Hence, for any given $m \in \mathbb{N}$ we have that $A^{\frac{1}{m}}$ and $B^{\frac{1}{m}}$ commute. Therefore $(AB)^{\frac{1}{m}}= A^{\frac{1}{m}}B^{\frac{1}{m}} \geq 0$. \end{proof}
\bigskip
It is shown in \cite{Inverse$M$-matrices1} that inverse $M$-matrices are closed under multiplication and addition by strictly positive diagonal matrices. It is natural to inquire if this property extends to every infinitely divisible matrix. In other words, if a nonnegative matrix $B$ is infinitely divisible if and only if $BD$ and $B+D$ are infinitely divisible for every strictly positive diagonal matrix $D$. It turns out that infinitely divisible matrices are not closed under positive diagonal multiplication. For example, consider the infinitely divisible matrix

\begin{equation} 
\begin{pmatrix} 2/5 & 2/5 & 1/5 \\ 0 & 1/2 & 1/2 \\ 0 & 0 & 1 \end{pmatrix}. 
\end{equation}

We can multiply this by a diagonal matrix to get,

\begin{equation} \label{counterexample}
\begin{pmatrix} 2/5 & 2/5 & 1/5 \\ 0 & 1/2 & 1/2  \\ 0 & 0 & \frac{1}{2} \end{pmatrix}. 
\end{equation} 

\noindent Computing the principal logarithm of matrix \eqref{counterexample} we see it is not infinitely divisible. We can however establish the following weaker results.\\

\begin{lem} \label{permutation lemma} 
 Let $B$ be a strongly infinitely divisible matrix and $L$ a strictly positive monomial matrix. Then $L^{-1} B L$ is strongly infinitely divisible. This implies $LB$ is strongly  infinitely divisible if and only if $BL$ is strongly infinitely divisible. 
  \end{lem}

\begin{proof}
Let $B=e^{-Q}$ for some $Z$ matrix $Q$. Clearly $L^{-1} B L=L^{-1} e^{-Q} L=e^{-L^{-1} Q L}$. Since $Q$ is a $Z$ matrix,  $L^{-1} Q L$ is a $Z$ matrix because $L^{-1}$ is non-negative and diagonal elements are only mapped to diagonal elements. Hence, $L^{-1} B L$ must be infinitely divisible. \\

Now assume that $LB$ is strongly infinitely divisible, then so is $L^{-1} L B L=BL$. A similar argument applies if we assume that $BL$ is strongly infinitely divisible. \end{proof}
With regards to the embedding problem, the case of interest is when $L$ is a permutation matrix. The above Lemma implies that if $P$ is an embeddable stochastic matrix, then so is $L P L^{-1} = L P L^{T}$. We can now prove Theorem \ref{dichotomy theorem}.\\

\begin{proof}[Proof of Theorem \ref{dichotomy theorem}] The matrix strongly infinitely divisible matrix $B$ is either irreducible or reducible, if it is irreducible, then it is known that there exists an $m \in \mathbb{N}$ such that $B^{m}$ is strictly positive. However this violates the invariance of zero patterns of Corollary \ref{zero pattern invariance} unless $B$ is strictly positive. This yields the first statement of the dichotomy \\

If the matrix $B$ is reducible then the decomposition $B=L U L^{T}$ as stated in Theorem \ref{dichotomy theorem} $ii)$ is a well known fact, see \cite{MatrixAnalysis1}. Furthermore, in this decomposition, the diagonal block matrices must be irreducible and thus by the same argument as in the preceding paragraph, this diagonal block matrices must be strictly positive.\\

We now show that the submatrices $U^{(n)}$, as defined in the statement of Theorem \ref{dichotomy theorem}, must be infinitely divisible. Without loss of generality assume $U$ is upper block triangular and that the eigenvalues are distinct. By Lemma \ref{permutation lemma}, $U$ is infinitely divisible. Consider the associated $Z$-matrix $-Q$ and the  submatrix of $-Q$, denoted $-Q^{(n)}$, obtained by deleting the first $n$ blocks from the top rows and left columns. Because $Q$ must be a polynomial of $U$, $-Q^{(n)}$ depends only on the entries in $U^{(n)}$. Therefore we conclude that that $-Q^{(M)}$ is a $Z$ matrix and $e^{-Q^{(n)}}=U^{(n)}$. \end{proof}

In light Theorem \ref{dichotomy theorem}, whenever dealing with strongly infinitely divisible nonnegative matrices, we may without loss of generality assume that it is strictly positive or upper block triangular. 

\bigskip

\section{The embedding problem for finite state stationary Markov chains}


The embedding problem for Markov chains has been a long standing problem in linear algebra and probability theory since it was first considered by Elfving \cite{Elfving}. It raises the question if a given discrete finite state Markov chain can be interpreted as having arisen from a continuous stationary Markov chain that has been observed at discrete intervals. Such Markov chain is called embeddable. This problem has found applications in a diverse number of fields, such as sociology \cite{Singer1976}, credit ratings \cite{Israel2001} and biology \cite{Nucleotide}. A Markov chain with stochastic Matrix $P$ is embeddable if and only if there exists an intensity matrix $R$ such that

\begin{equation} P=e^{R} \end{equation}

The reader is directed to Singer and Spilerman \cite{Singer1976} for the definition of an intensity matrix and a wide variety of of examples illustrating the depth of this problem. Kingman \cite{Kingman1962}  showed that a Markov chain was embeddable if and only if it was nonsingular and had stochastic matrix roots of arbitrary order. We thus recognize the embeddable matrices as a special case of strongly infinitely divisible matrices.\\

For convenience and clarity we note what our key results entail for the embedding problem for stochastic matrices. Before this, however, there are a few things to verify. The following was proved recently by EB Davies \cite{EBDavies} and can be proved in a similar way to Lemma \ref{infdivdense}.\\

\begin{thm}
 (EB Davies, 2010) The set of matrices with distinct eigenvalues dense in the set of embeddable stochastic matrices.
\end{thm}

In light of the above theorem, we now  realize the implication of Lemma \ref{infdivdense} on the embedding problem to deduce what was and can show what was proved by Kingman \cite{Kingman1962} without the additional assumption that the roots are also stochastic. \footnote{A stochastic matrix may have nonnegative nonstochastic roots, an example is given in \cite{Higham2011}} \\

\begin{thm} Assume a stochastic matrix $P$ is strongly infinitely divisible, then $P$ is embeddable. \end{thm}

\begin{proof}
It suffices to consider  a stochastic matrix $P$ with distinct eigenvalues. \\
 By Theorem \ref{Z-matrix theorem}, if $P$ has nonnegative roots for all $n$. Then there is a $Z$ matrix $Q$ such that

\begin{equation} P=e^{-Q}. \end{equation}

 Let $u$ be the vector of length $N$, all of whose entries are $1$. Since $P$ has distinct eigenvalues,  $Q=\log(P),$ is a polynomial of $P$. It follows that $u$ is a eigenvector of $Q$, with eigenvalue $0$. However because $Q$ is a $Z$ matrix, this implies that $-Q$ must in fact be a intensity matrix and hence $P$ is embeddable. \end{proof}
 
 \bigskip
 
 \indent It is useful to note that, in the case of stochastic matrices, these $M$-matrices must be of a specific form. If a stochastic matrix $P$ is the inverse of an $M$ matrix then we have that $P^{-1}=sI-K$ for some $s > \rho(K)$ and $K$ is nonnegative. However, we know that $u$, the vector consisting of ones as defined above, must be an eigenvector for $K$. Let $\lambda = \rho(K)$, so that  $\lambda=s-1$. Defining $H=\frac{1}{\lambda} K$, where $H$ is now stochastic, we have

 \begin{equation} 
 P^{-1}=s \Big (I-\frac{s-1}{s} H \Big ). 
 \end{equation}

 \noindent Thus when an embeddable stochastic matrix $P$ is the power of an inverse $M$ matrix, it must be of the form

  \begin{equation} P=(1-\epsilon)^{m} (I-\epsilon H)^{-m}. \end{equation}

  \noindent Where $\epsilon=\frac{s-1}{s}$.  We can thus classify the set of embeddable stochastic matrices: those stochastic matrices which can be infinitesimally perturbed to be in the form $(1-\epsilon)^{m} (I-\epsilon H)^{-m}$. More formally, we state the following result. \\

 \begin{cor} A stochastic matrix is $P$ is embeddable if and only if it is nonsingular and in the closure of the set:  
 
 \begin{equation} \lbrace P \: : \: P=(1-\epsilon)^{m} (I-\epsilon H)^{-m}, \: 0 \leq \epsilon < 1, \: m \in \mathbb{N}, \: H \:\: \textnormal{stochastic} \rbrace. \end{equation}
  \end{cor}
 \bigskip \bigskip

Bounds on the eigenvalues for the intensity matrices of embeddable stochastic matrices have been developed in some length. Notably, Runnenbergs' condition \cite{Runnenberg1962} which states that the eigenvalues of a values of an $n \times n$ intensity matrix must be an element of the set

\begin{equation}
 \big \lbrace z \: : \: \pi \Big (\frac{1}{2}+\frac{1}{n} \Big ) \leq \textnormal{arg} (z) \leq \Big ( \frac{3}{2}-\frac{1}{n} \Big ) \big \rbrace.
\end{equation}

\noindent The utility of this  result however, diminishes rapidly in higher dimensions. More in the spirit of this analysis is the related bound proven in \cite{Israel2001}:

\begin{equation}  \label{Israelbound}
| \textnormal{Im} (\log (\lambda)) | \leq - \log(\det (P)). 
\end{equation}

 \noindent If we apply our bound derived at the end of Section $1$ to the case of stochastic matrices, we arrive at  \eqref{Israelbound}.\\
 
  We also know that the diagonal elements in any intensity matrix are always nonpositive and the rows sum to $0$. Hence  using Gershgorin's disc theorem, the imaginary part of any eigenvalue of an intensity matrix is nonpositive; therefore 
 
 \begin{equation} \label{mybound}
  0 \geq \textnormal{Im} \log (\lambda) \geq \log(\det(P))
   \end{equation}
 
 \noindent Inequality \eqref{mybound} dramatically simplifies the procedure for determining whether a stochastic matrix is embeddable. For practical purposes we can usually, without lost of generality, restrict ourselves to the case of distinct eigenvalues, as this may always be obtained after a infinitesimal perturbation by Lemma \ref{infdivdense}. In this case one needs only to check branches of the logarithm with imaginary part in the domain above. For example, if $P$ is a $5 \times 5 $ matrix, and the determinant of $P$ is small, say $0.00001$, we need to check only $16$ cases.
 \bigskip

Theorem \ref{dichotomy theorem} also has a probabilistic interpretation: this result implies that there are in fact only two types of finite state, stationary, continuous Markov chains. One type corresponds to a process whereby, from any state, it may, with positive probability, reach any other state in any given time interval. This type corresponds to Theorem 3 (i). The other type is when there is a hierarchy of systems described by  some sequence of square upper triangular block matrices $P^{n}$, $n \leq M$, each modeling a continuous Markov chain in its own right. 
The practical application of 
this result is that one can can determine if a stochastic matrix is embeddable by checking if the stochastic matrices defined by submatrices $P^{n}$ are embeddable. I.e our result introduces a new necessary condition.\\

 \section{Acknowledgments}

This research was supported by a Victoria University of Wellington Summer Scholarship and a Monbukagakusho scholarship administered by the Japanese Ministry of Education, Culture, Sports, Science and Technology. I would also like to thank Matt Visser for his helpful insight and guidance, in particular with regards to the counterexamples \eqref{countexam1} and \eqref{countexam2}. \.  Finally, I would like to thank a patient referee for the helpful comments.

\pagebreak

\nocite{Higham2011} 
\nocite{InverseM-matrices2}
 \textbf{References}
\bigskip
\bibliography{Linearalgebra}

\begin{thebibliography}{10}

\bibitem{Plemmons1977}
R.~Plemmons, ``{M-matrix characterizations. I-nonsingular $M$-matrices},'' {\em
  Linear Algebra and its Applications}, vol.~no.18, no.~2, pp.~pp.175 -- 188,
  1977.

\bibitem{MatrixAnalysis2}
R.~A. Horn and C.~R. Johnson, {\em Topics in Matrix Analysis}.
\newblock Cambridge University Press, 1994.

\bibitem{Inverse$M$-matrices1}
C.~R. Johnson, ``{Inverse $M$-matrices},'' {\em Linear Algebra and its
  Applications}, vol.~no.47, pp.~pp.195--216, 1982.

\bibitem{InverseM-matrices2}
C.~R. Johnson and R.~L. Smith, ``{Inverse M-matrices \{II\}},'' {\em Linear
  Algebra and its Applications}, vol.~435, no.~5, pp.~pp.953 -- 983, 2010.

\bibitem{Elfving}
E.~G, ``{Zur Theorie der Markoffschen},'' {\em Acta Soc. Sci Fennicae n. Ser A
  2}, vol.~no.8, 1937.

\bibitem{Kingman1962}
J.~Kingman, ``{The imbedding problem for finite Markov chains},'' {\em
  {Zeitschrift für Wahrscheinlichkeitstheorie und Verwandte Gebiete}},
  vol.~no.1, pp.~pp.14--24, 1962.

\bibitem{RAHorn1969}
R.~A.Horn., ``{The theory of infinitely divisible matrices and kernels},'' {\em
  Trans.Amer. Math. Soc.}, vol.~136, pp.~pp.269--286, 1969.

\bibitem{MatrixFunction}
N.~J. Higham, {\em Functions of Matrices: Theory and Computation}.
\newblock SIAM, 2008.

\bibitem{MatrixAnalysis1}
R.~A. Horn and C.~R. Johnson, {\em Matrix Analysis}.
\newblock Cambridge University Press, 1990.

\bibitem{Singer1976}
B.~Singer and S.~Spilerman, ``{The Representation of Social Processes by Markov
  Models},'' {\em American Journal of Sociology}, vol.~no.82, no.~1,
  pp.~pp.1--54, 1976.

\bibitem{Denseeigenvalues}
D.~J. Hartfiel, ``{Dense Sets of Diagonalizable Matrices},'' {\em Proceedings
  of the American Mathematical Society}, vol.~no.123, no.~6, pp.~pp.
  1669--1672, 1995.

\bibitem{Israel2001}
R.~B. I. J.~S. Roosenthal{,} and J.~Z. Wei, ``{Finding Generators for Markov
  Chains via Empirical Transition Matrices, with Applications to Credit
  Ratings},'' {\em Mathematical Finance}, vol.~no.11, pp.~pp.245--265, 2001.

\bibitem{$M$-MatrixFunctions}
M.~Fiedler and H.~Schnieder, ``{Analytic Functions of $M$-Matrices and
  Generalisations. },'' {\em Linear and multilinear algebra}, vol.~no.13,
  pp.~pp.185--201, 1983.

\bibitem{Nucleotide}
V.~K. Y. V. P. A.~S. Y and H.~GA, ``{The Embedding Problem for Markov Models of
  Nucleotide Substitution},'' {\em PloS ONE}, 2013.

\bibitem{EBDavies}
E.~B. Davies, ``{Embeddable Markov Matrices},'' {\em Electronic Journal of
  Probability}, vol.~no.15, no.~47, pp.~pp.1474--1486, 2010.

\bibitem{Higham2011}
N.~J. Higham and L.~Lin, ``{On pth roots of stochastic matrices},'' {\em Linear
  Algebra and its Applications}, vol.~no.435, pp.~pp.448 -- 463, 2011.

\bibitem{Runnenberg1962}
R.~J, ``{On Elfving’s problem of imbedding a time-discrete Markov chain in a
  time continuous one for finitely many states},'' {\em Indagationes
  Mathematicae (Proceedings)}, 1962.

\end{thebibliography}


\begin{thebibliography}{1}

\bibitem{VANBRUNT2018163}
A.~Van-Brunt, ``Infinitely divisible nonnegative matrices, m-matrices, and the
  embedding problem for finite state stationary markov chains,'' {\em Linear
  Algebra and its Applications}, vol.~541, pp.~163 -- 176, 2018.

\bibitem{MatrixFunction}
N.~J. Higham, {\em Functions of Matrices: Theory and Computation}.
\newblock SIAM, 2008.

\bibitem{MatrixAnalysis1}
R.~Horn and C.~Johnson, {\em Matrix Analysis}.
\newblock Cambridge University Press, 2013.

\bibitem{Higham2011}
N.~J. Higham and L.~Lin, ``{On pth roots of stochastic matrices},'' {\em Linear
  Algebra and its Applications}, vol.~no.435, pp.~448 -- 463, 2011.

\bibitem{InverseM-matrices2}
C.~R. Johnson and R.~L. Smith, ``{Inverse M-matrices \{II\}},'' {\em Linear
  Algebra and its Applications}, vol.~435, no.~5, pp.~953 -- 983, 2010.

\end{thebibliography}

\end{document}